\documentclass[conference]{IEEEtran}
\IEEEoverridecommandlockouts
\usepackage{cite}
\usepackage{algorithmic}
\usepackage{textcomp}
\def\BibTeX{{\rm B\kern-.05em{\sc i\kern-.025em b}\kern-.08em
		T\kern-.1667em\lower.7ex\hbox{E}\kern-.125emX}}

\usepackage{amssymb, amsmath,amsfonts}
\usepackage{amsbsy}

\usepackage[english]{babel}
\usepackage{accents}
\usepackage{dsfont}
\usepackage{booktabs}
\usepackage{hyperref}
\hypersetup{
	colorlinks,
	linkcolor={red!50!black},
	citecolor={blue!50!black},
	urlcolor={blue!80!black}
}

\usepackage{etoolbox}
\apptocmd{\sloppy}{\hbadness 10000\relax}{}{}

\usepackage{enumerate}
\usepackage[shortlabels]{enumitem}
\usepackage{ifthen}
\usepackage{float}
\floatstyle{ruled}
\newfloat{algorithm}{htb}{alg}
\floatname{algorithm}{Algorithm}
\newcounter{algo@row}
\newcounter{algo@rowindent}
\newcommand{\algofont}[1]{\textbf{#1}}
\newcommand{\algonumbersize}[1]{\scriptsize{#1}}
\newcommand{\algopreitem}[1][\arabic{algo@row}]{\texttt{\algonumbersize{#1}}}
\newcommand{\algoitemskip}{\hspace{\value{algo@rowindent}cc}}
\newenvironment{algo}{\vskip.3em\small%
	\begin{list}{\algopreitem\texttt{\algonumbersize{:}}}{%
			\usecounter{algo@row}%
			\setcounter{algo@rowindent}{0}%
			\setlength{\itemindent}{2em}%
			\setlength{\labelwidth}{2em}
			\setlength{\parsep}{0cm}%
		}%
	}{
	\end{list}
}
\newcommand{\algonewnestedopen}[2]{
	\newcommand{#1}[1][]{%
		\ifthenelse{\equal{##1}{}}{\item}{\item[{\algopreitem[##1]}]}
		\algoitemskip\algofont{#2}%
		\addtocounter{algo@rowindent}{1}%
		\ignorespaces
	}
}
\newcommand{\algonewnestedaux}[2]{
	\newcommand{#1}[1][]{
		\addtocounter{algo@rowindent}{-1}
		\ifthenelse{\equal{##1}{}}{\item}{\item[{\algopreitem[##1]}]}
		\algoitemskip\algofont{#2}%
		\addtocounter{algo@rowindent}{+1}%
		\ignorespaces
	}
}
\newcommand{\algonewnestedclose}[2]{
	\newcommand{#1}[1][]{
		\addtocounter{algo@rowindent}{-1}
		\ifthenelse{\equal{##1}{}}{\item}{\item[{\algopreitem[##1]}]}
		\algoitemskip\algofont{#2}%
		\ignorespaces
	}
}
\newcommand{\algonewcommand}[2]{
	\newcommand{#1}[1][default]{
		\ifthenelse{\equal{##1}{default}}{\item}{\item[{\algopreitem[##1]}]}%
		\algoitemskip\algofont{#2}%
		\ignorespaces
	}%
}
\newcommand{\algonewkeyword}[2]{\newcommand{#1}{\algofont{#2}}}
\algonewcommand{\STATE}{\ignorespaces}
\algonewcommand{\INPUT}{Input: }
\algonewcommand{\pINPUT}{\phantom{Input: }}
\algonewcommand{\COMPUTE}{Compute: }
\algonewcommand{\OUTPUT}{Output: }
\algonewcommand{\pOUTPUT}{\phantom{Output: }}
\algonewnestedopen{\IF}{if }
\algonewnestedaux{\ELSEIF}{else if }
\algonewnestedaux{\ELSE}{else }
\algonewnestedclose{\ENDIF}{end if }
\algonewnestedopen{\FOR}{for }
\algonewnestedclose{\ENDFOR}{end for }
\algonewnestedopen{\PARFOR}{parfor }
\algonewnestedclose{\ENDPARFOR}{end parfor }
\algonewnestedopen{\WHILE}{while }
\algonewnestedclose{\ENDWHILE}{end while }
\algonewcommand{\BREAK}{break}%
\algonewkeyword{\For}{for }%
\algonewkeyword{\To}{to }%
\algonewkeyword{\Do}{do }%
\algonewkeyword{\If}{if }%
\algonewkeyword{\Then}{then }%
\algonewkeyword{\Else}{else }%
\algonewkeyword{\End}{end }%
\algonewkeyword{\AND}{and }%
\algonewkeyword{\True}{true }%
\algonewkeyword{\False}{false }%
\algonewkeyword{\Call}{call }%
\algonewkeyword{\irbleigs}{irbleigs }%
\algonewkeyword{\tridiag}{tridiag}%
\algonewkeyword{\reorth}{reorth}%

\usepackage{multirow}

\usepackage{graphicx}
\usepackage[dvipsnames,x11names]{xcolor}

\usepackage{todonotes}

\usepackage{amsbsy}

\newcommand{\K}{{\cal K}}
\newcommand{\R}{{\cal R}}

\newcommand{\N}{\mathbb{N}}
\newcommand{\bsf}{{\boldsymbol f}}
\newcommand{\bsg}{{\boldsymbol g}}
\newcommand{\MSRE}{{\boldsymbol{\mathcal{E}}}}

\DeclareMathSymbol{\shortminus}{\mathbin}{AMSa}{"39}

\DeclareMathOperator*{\argmin}{argmin}

\makeatletter
\newcommand{\sbullet}{%
	\hbox{\fontfamily{lmr}\fontsize{.6\dimexpr(\f@size pt)}{0}\selectfont\textbullet}}

\makeatother

\newcommand{\Kd}{{\mathcal{K}^{(n)}}}

\date{}

\begin{document}
\title{Graph approximation and generalized Tikhonov regularization for signal deblurring
	\thanks{Supported by INdAM-GNAMPA and INdAM-GNCS.}}

\author{\IEEEauthorblockN{Davide Bianchi}
	\IEEEauthorblockA{\textit{dept. of Science and High Technology} \\
		\textit{University of Insubria}\\
		Como, Italy \\
		d.bianchi9@uninsubria.it}
	\and
	\IEEEauthorblockN{Marco Donatelli}
	\IEEEauthorblockA{\textit{dept. of Science and High Technology} \\
		\textit{University of Insubria}\\
		Como, Italy \\
		marco.donatelli@uninsubria.it}
	}

\maketitle              
\begin{abstract}
Given a compact linear operator $\K$, the (pseudo) inverse $\K^\dagger$ is usually substituted by a family of regularizing operators $\R_\alpha$ which depends on $\K$ itself. Naturally, in the actual computation we are forced to approximate the true continuous operator $\K$ with a discrete operator $\K^{(n)}$ characterized by a finesses discretization parameter $n$, and obtaining then a discretized family of regularizing operators $\R_\alpha^{(n)}$. In general, the numerical scheme applied to discretize $\K$ does not preserve, asymptotically, the full spectrum of $\K$. In the context of a generalized Tikhonov-type regularization, we show that a graph-based approximation scheme that guarantees, asymptotically, a zero maximum relative spectral error can significantly improve the approximated solutions given by $\R_\alpha^{(n)}$. This approach is combined with a graph based regularization technique with respect to the penalty term.
\end{abstract}

\section{Introduction}
The aim of this work is to provide a first glimpse about the applications of a full graph-based approximation approach  to inverse problems regularization. Since the theory involved can be rather technical and very vast, to keep this manuscript short and almost self-contained, we will mainly concentrate on numerical examples in the one dimensional case. 

That said, the model equation we consider is
\begin{equation}\label{model_problem}
 \K [f] = g,
\end{equation}
where $\K : L^2\left([0,1]\right) \to  L^2\left([0,1]\right)$ is a compact linear operator acting on the Hilbert space of square integrable functions over $[0,1]$. In particular, $\K$ will be the Green operator of a self-adjoint second-order differential operator $\mathcal{L}$ with formal equation
\begin{equation*}
\mathcal{L}[g](x):= -g''(x) + q(x)g(x) \qquad \forall x\in (0,1).
\end{equation*} 
Indicating with $\K^\dagger$ the generalized Moore-Penrose inverse of $\K$, the solution for the model problem \eqref{model_problem} of minimal norm reads $f^\dagger = \K^\dagger g$. It is known that \eqref{model_problem} is an ill-conditioned problem, that is, even small errors in the observed data $g$ will greatly affect the reconstructed solution. 

We are interested to approximate the solution $f^\dag$ when only a noisy approximation $g^\epsilon:= g +\eta$ is available, with
\begin{equation*}\label{eq:delta}
\|g^\epsilon - g \|=\|\eta \| = \epsilon,
\end{equation*}
and where $\epsilon$ is called the \emph{noise level}.  Since $\K^\dag g^\epsilon$ is not a good
approximation of $f^\dag$ due to the ill-conditioning of $\K$, it is commonly chosen to approximate $f^\dag$ with
$f_\alpha ^\epsilon := \R_\alpha [g^\epsilon]$,
where $\{\R_\alpha\}_{\alpha \in (0,+\infty)}$ is a family of continuous operators depending on a parameter
$\alpha$ such that $\R_\alpha \to \K^\dagger$ pointwise as $\alpha=\alpha(\epsilon) \to 0$.
A classical example is the generalized Tikhonov regularization method defined by
\begin{subequations}
\begin{equation}\label{def:gen_tik1}
\mathcal{R}_{\alpha} : L^2([0,1])\to \textnormal{dom}\left(\mathcal{A}\right)\subseteq L^2([0,1]) \quad \mbox{such that }
\end{equation}
\begin{equation}\label{def:gen_tik2}
\mathcal{R}_{\alpha} [g^\epsilon]:= \textnormal{argmin}_{f \in \textnormal{dom}\left(\mathcal{A}\right)} \left\{ \left\| \mathcal{K}[f] - g^\epsilon\right\|^2 + \alpha \left\|\mathcal{A} [f] \right\|^2\right\},
\end{equation}  
\end{subequations}	
where $\mathcal{A} :  \textnormal{dom}\left(\mathcal{A}\right)\subseteq L^2([0,1]) \to L^2([0,1])$ is a closed and densely defined linear operator such that
\begin{equation*}
\ker\left(\mathcal{A}\right)<\infty, \qquad \ker\left(\mathcal{A}\right)\cap \ker\left(\mathcal{K}\right)= \{0\}.
\end{equation*}
$\mathcal{R}_\alpha$ is called \emph{generalized Tikhonov regularization operator} and we will write $f_\alpha^\epsilon := \mathcal{R}_{\alpha} \left[g^\epsilon\right]$ for the regularized solution corresponding to the given data $g^\epsilon$. For a detailed account over this topic we refer to \cite{G,EHN}. For $\mathcal{A}=\mathcal{I}$, where $\mathcal{I}$ is the identity map, then we recover the standard Tikhonov regularization. The operator $\mathcal{A}$ is typically introduced to force the regularized solution $f_\alpha^\epsilon$ to live in $\textnormal{dom}\left(\mathcal{A}\right)$, whenever we have a-priori informations on particular features of the true solution $f^\dagger$.

To better understand the role of $\mathcal{A}$ and to help us to simplify the computations, let us make the following, momentarily, assumptions:  $\mathcal{A}$ is self-adjoint (with purely discrete spectrum) and it shares the same eigenbase with $\K$. Indicating with $\{\lambda_m; u_m, v_m\}_{m\in \N}, \{\mu_m; u_m, v_m\}_{m\in \N}$ the spectral decomposition of $\K$ and $\mathcal{A}$, respectively, we can then express the regularized solution $f_\alpha^\epsilon$ as
$$
f_\alpha^\epsilon = \sum_{m\in \N} \frac{\lambda_m}{\lambda_m^2 + \alpha\mu_m^2 } \langle g^\epsilon, u_m \rangle u_m.
$$  
If $\ker(\mathcal{A})\neq \{0\}$, then we can write
\begin{align*}
f_\alpha^\epsilon = &\sum_{\substack{u_m \in \ker(\mathcal{A})}}\lambda_m^{-1}\langle g^\epsilon, u_m \rangle u_m \\
&+ \sum_{\substack{u_m \notin \ker(\mathcal{A})}} \frac{\lambda_m}{\lambda_m^2 + \alpha\mu_m^2 } \langle g^\epsilon, u_m \rangle u_m.
\end{align*}
The above equation tells that the regularized solution $f_\alpha^\epsilon$ is made of the sum of two parts: the first part is the projection of the observed data $g^\epsilon$ into $\ker(\mathcal{A})$ and no regularization takes action since $\alpha$ does not play any role in it, and the second part is the remainder of the summation where instead the regularization operator takes action. Writing $F_\alpha (t) :=\frac{t^2}{t^2 + \alpha}$, and
\begin{align*}
& \textbf{err}_{app,\perp} := -\sum_{u_m \notin \ker(\mathcal{A})} \left(1 - F_{\alpha} \left(\frac{\lambda_m}{\mu_m}\right) \right)\lambda_{m} ^{-1} \langle g, u_m \rangle u_m,\\
&\textbf{err}_{noise,\perp} := \sum_{u_m \notin \ker(\mathcal{A})}F_{\alpha} \left(\frac{\lambda_m}{\mu_m}\right) \lambda_{m} ^{-1} \langle \eta, u_m \rangle u_m,
\end{align*}
then 
\begin{align*}
\|f_\alpha^\epsilon - f^\dagger\| &= \left\|f_\alpha^\epsilon - f^\dagger\right\|_{\ker(\mathcal{A})} + \left\|f_\alpha^\epsilon - f^\dagger\right\|_{\ker(\mathcal{A})^\perp}\\
&= \left\|\sum_{u_m \in \ker(\mathcal{A})}\lambda_{m} ^{-1}\left(\langle g^\epsilon, u_m \rangle u_m  -  \langle g, u_m \rangle u_m \right)\right\|\\
&\,  +  \left\| \textbf{err}_{app,\perp} + \textbf{err}_{noise,\perp}\right \|\\
&= \left\|\sum_{u_m \in \ker(\mathcal{A})}\lambda_{m} ^{-1}\langle \eta, u_m \rangle u_m   \right\|\\
&\, +  \| \textbf{err}_{app,\perp} + \textbf{err}_{noise,\perp} \|\\
&\leq c\epsilon  + \| \textbf{err}_{app,\perp}  + \textbf{err}_{noise,\perp} \|,
\end{align*}
where $c= \max\{\lambda_{m} ^{-1} \, : \, u_m \in \ker(\mathcal{A}) \}$. So, if $f^\dagger$ belongs entirely to $\ker(\mathcal{A})$, the best possible strategy would be to just project $f_\alpha^\epsilon$ onto the kernel of $\mathcal{A}$, such that to delete the error coming from $\textbf{err}_{app,\perp}$ and $\textbf{err}_{noise,\perp}$. Morally speaking, when choosing the operator $\mathcal{A}$ we should look that as much features as possible of $f^\dagger$ belong to  $\ker\left(\mathcal{A}\right)$. For example, a typical choice for $\mathcal{A}$ in imaging is the Laplace operator with Dirichlet or Neumann boundary conditions: the first one is chosen when dealing with astronomical images that have a black background, which correspond to a zero numeric value at the boundary of the field of view, in the grayscale representation; the second one is chosen because images have several areas of homogeneous, constant color. For a reference, see \cite{AH}. Clearly, in the real applications it is not usually granted that $\K$ and $\mathcal{A}$ commute, and in some works it was proposed to use $\mathcal{A}=F\left(\K\K^*\right)$, where $F$ is a suitable function on the spectrum of $\K\K^*$ that mimic the spectral distribution of the Laplace operator, see \cite{HS,LNEM,BD}.

In the actual computations we can not make use of the continuous operator $\K$ and $\mathcal{A}$, but only of discrete approximations $\Kd$ and $\mathcal{A}^{(n)}$, respectively, that are obtained by numerical discretization schemes. As a consequence, the ideal family of regularization operators $\left\{\R_\alpha\right\}_\alpha$ is substituted by a family of discrete operators $\left\{\R_\alpha^{(n)}\right\}_\alpha$, such that 
 \begin{equation}\label{def:discrete_generalized_Tikhonov}
 \R_\alpha^{(n)} [\bsg_n^{\epsilon}] := \argmin_{\bsf_n\in L_n} \left\{ \|\Kd[\bsf_{n}] - \bsg_n^{\epsilon}\|^2 + \alpha \|\mathcal{A}^{(n)}[\bsf_{n}]\|^2  \right\},
 \end{equation}
where $L_n$ is a finite dimensional space that approximates $L^2([0,1])$. Clearly, the discrete approximation can introduce extra round-off errors, part of which comes from a bad approximation of the full spectrum of $\K$. The idea is that if we can have a better control on the spectral relative error, we can obtain a better regularized solution.

This paper is organized as follows: in Section \ref{sec:graph_def} we report some notations from graph theory. In Section \ref{sec:graph_approx} we show a graph based approximation method to obtain a discretized operator such that it preserves uniformly the original spectrum of $\K$. We will see that the key point is to use the Fourier coefficients of the inverse Weyl distribution function associated to $\K$ (or, equivalently, to $\mathcal{L}$) as weights of a transformed-path graph Laplacian. In Section \ref{sec:graph_penalty}, we build the operator $\mathcal{A}$ in the penalty term as the graph Laplacian associated to a given graph $G$. The main idea is to encode in $\mathcal{A}$ informations from the observed data $g^\epsilon$ that can help to better define the space where to force into the regularized solution. Finally, in Section \ref{sec:numerical_examples} we provide some numerical examples.

\section{Graph setting}\label{sec:graph_def}
Given a countable set of nodes $X=\left\{x_i : i\in I\right\}$, a \textit{graph} $G$ on $X$ is a pair $(w,\kappa)$ such that
\begin{itemize}
	\item $w: X\times X \to [0,\infty)$ is a nonnegative, symmetric function with zero diagonal; 
	\item $\kappa : X \to [c,\infty)$ is a lower-bounded function, $c>-\infty$.
\end{itemize}
In a more compact notation, we write $G=(X,w,\kappa)$, and whenever $\kappa\equiv 0$ we will just write $G=(X,w)$. $k$ is called \emph{potential} (or \emph{killing}) term. For a detailed modern introduction to graph theory, see \cite{KLW_book}. Every unordered pair of nodes $e=\{x_i,x_j\}$ such that $w(x_i,x_j)> 0$ is called \emph{edge} incident to $x_i$ and to $x_j$, and the collection $E$ of all the edges is uniquely determined by $w$. The non-zero values $w(x_i,x_j)$ are called \emph{weights} associated with the edges $\{x_i,x_j\}$, and $w$ is the \emph{edge-weight} function. Two nodes $x_i,x_j$ are said to be \textit{neighbors} (or \emph{connected}) in $G$ if $\{x_i,x_j\}$ is an edge and we write $x_i\sim x_j$. In this work, we will assume that $X$ is endowed with the counting measure and that every node $x_i$ is connected at most with a finite number of nodes $x_j$. The sum of all the weights incident to a node $x_i$ is called the \emph{degree} of $x_i$. In the case that $w(x_i,x_j)\in \{0,1\}$, then $G$ is said to be \emph{unweighted}. Given a real valued function $f$ on $X$, the (formal) \emph{graph Laplacian} $\Delta$ of $G$ applied to $f$ is defined via
\begin{equation}\label{formal_laplacian}
\Delta[f](x_i):= \sum_{\substack{x_j \in X}} w(x_i,x_j)\left(f(x_i) - f(x_j)\right) + \kappa(x_i)f(x_i).
\end{equation}
For a \emph{connected} graph, that is, a graph where for every pair $\{x_i, x_j\}$ there exists a sequence of connected nodes such that $x_i=x_{i_0}\sim \cdots \sim x_{i_k}=x_j$, the \emph{combinatorial} distance function $d : X\times X \to [0,\infty)$ reads as
\begin{equation*}
d(x_i,x_j):= \min\{k \in \N \, : \,  x_i= x_{i_0}\sim x_{i_1}\sim \cdots \sim x_{i_k}=x_j \},
\end{equation*}
and $d_\infty := \max \{ d(x_i,x_j)\, : \, x_i,\, x_j \in X\}$ is the \emph{combinatorial diameter} of $G$. Given an unweighted graph $G=(X,w)$ and a proper subset $X_0\subset X$, we call \emph{Dirichlet $m$-path graph Laplacian}, associated to $X_0$ and $G$, the operator $\Delta_{k,\textnormal{dir}}$ defined by
\begin{equation}\label{def:k-path-Laplacian}
\Delta_{m,\textnormal{dir}}[f](x_i) := \sum_{\substack{x_j \in X_0:\\d(x_i,x_j)=m}} \left(f(x_i) - f(x_j)\right) + \kappa_{m,\textnormal{dir}}(x_i)f(x_i),
\end{equation}
where 
\begin{equation*}
\kappa_{m,\textnormal{dir}}(x_i):= \sum_{x_j \in X\setminus X_0} \mathds{1}_{\{k\}}\left(d(x_i,x_j)\right)
\end{equation*}
and
\begin{equation*}
\mathds{1}_{\{k\}}(t):= \begin{cases}
1 & \mbox{if } t=k,\\
0 & \mbox{otherwise}.
\end{cases}
\end{equation*}
We call $\kappa_{m,\textnormal{dir}}(\cdot)$ the \emph{Dirichlet boundary potential}. Finally, given a function $\phi : \N \to \mathbb{R}$, we define the \emph{Dirichlet transformed-path graph Laplacian} $\Delta_{\infty,\textnormal{dir}}$ as
\begin{equation}\label{def:transformed-path-Laplacian}
\Delta_{\infty,\textnormal{dir}} := \sum_{m=1}^{d_\infty} \phi(m)\Delta_{m,\textnormal{dir}}.
\end{equation}
About \eqref{def:k-path-Laplacian} and \eqref{def:transformed-path-Laplacian}, they are a slight generalization of the original definitions of $m$-path graph Laplacian and the transformed-path graph Laplacian, respectively, which can be found in \cite{Estrada1,Estrada2}. Those generalizations take into account the potential term $\kappa_{m,\textnormal{dir}}$ that counts the edge deficiency of a node in $X_0$ with respect to its degree value as a node in $X$, the node-set of the original graph $G$, see \cite{KL} and \cite[Section 7.1]{ABS} for an example of graph discretization of a PDE by Dirichlet subgraphs.

\section{Graph approximation of $\K$}\label{sec:graph_approx}
Making a simplification, when discretizing problem \eqref{model_problem} we pass from the continuous operator $\K$ to a (matrix) discrete operator $\K^{(n)}$ which acts on a finite dimensional space $L_n$, isomorph to $\mathbb{R}^n$, such that $L_n \subset L_{n+1} \subset \ldots \subset L^2([0,1])$ and $\overline{\bigcup_{n=1}^\infty L_n} = L^2([0,1])$. For a detailed account about regularization methods by projection, we point the interested reader to \cite[Section 5.2]{EHN}, and all the references therein.

The discrete operator $\K^{(n)}$ has to be consistent with $\K$ and, therefore, at least it has to satisfy  
\begin{enumerate}[label={\upshape(\bfseries P\arabic*)},wide = 0pt, leftmargin = 3em]
	\item\label{item:local_convergence} $\| \K^{(n)} [\bsf_n] - \K[f] \| \to 0$ as $n\to \infty$, where $\bsf_n$ is the projection of $f$ into $L_n$ and $\|\cdot\|$, in this case, can be both the $L^2$ norm or the sup norm;
	\item\label{item:local_eig_estimate} $\left|\frac{\lambda^{(n)}_m}{\lambda_m} -1\right|\to 0$ as $n\to \infty$, for every \emph{fixed} $m \in \N$, where $\lambda_m ^{(n)}$ and $\lambda_m$ are the eigenvalues of the discretized operator $\Kd$ and of the continuous operator $\K$, respectively. 
\end{enumerate}
It often happens that estimate \ref{item:local_eig_estimate} does not work on all the spectrum, that is, if $m=m(n)$ is \emph{not} fixed then \ref{item:local_eig_estimate} is not satisfied. We call \emph{local spectral relative error} (LSRE) the numerical quantity in \ref{item:local_eig_estimate}. Writing 
\begin{equation}\label{def:MSRE}
\MSRE:= \limsup_{n\to \infty} \max_{m=1,\ldots,n}\left\{\left|\frac{\lambda^{(n)}_m}{\lambda_m} -1\right|\right\},
\end{equation}
the preceding remark translates into saying that, in general, $\MSRE>0$. Clearly, this introduces more errors in the regularization process. We call $\MSRE$ the \emph{maximum spectral relative error} (MSRE), see \cite{Bianchi} for more details. We want to show that, if we can guarantee 
\begin{enumerate}[label={\upshape(\bfseries P2')},wide = 0pt, leftmargin = 3em]
\item \label{item:uniform_eig_estimate} $\MSRE=0$,
\end{enumerate}
then we obtain a significant improvement in the approximation of the regularized solution $\bsf_n^\epsilon$. For simplicity, suppose that $\K$ is a compact linear operator such that
\begin{equation}\label{eq:Green_operator}
\K[f](x):= \int_0^1 h(x,y)f(y)\, dy 
\end{equation}
whose kernel $h$ is the Green function of the following second-order differential operator $\mathcal{L} : \textnormal{dom}\left(\mathcal{L}\right) \subset L^2([0,1]) \to L^2([0,1])$, 
\begin{equation}\label{eq:diff_operator}
\begin{cases}
\textnormal{dom}\left(\mathcal{L}\right):= H^1_0([0,1]),\\
\mathcal{L}[g](x):= -g''(x) + q(x)g(x), 
\end{cases}
\end{equation}  
where we suppose $q$ to be bounded, and $H^1_0$ is the usual closure of the set of compactly supported smooth functions $C_c^\infty(0,1)$ with respect to the Sobolev space $H^1$. In other words, $\mathcal{L}$ is a Schr\"{o}dinger-type operator with Dirichlet \emph{boundary conditions} (BCs). We suggest the following readings about this topic, \cite{Davies1,EM}. We can think of $\K$ as the pseudo-inverse of $\mathcal{L}$, i.e., $\K = \mathcal{L}^\dagger$. Therefore, if we find a way to discretize $\mathcal{L}$ such that the MSRE associated to the discretization of $\mathcal{L}$ is zero, i.e., such that $\MSRE=\MSRE(\mathcal{L})=0$, then we will have that the MSRE associated to the discretization of $\K$ will be zero, i.e., $\MSRE=\MSRE(\mathcal{K})=0$, by the substitution $\Kd = \mathcal{L}^{(n),\dagger}$.

For this regard, let $G=(X,w)$ be such that
\begin{itemize}
	\item $X=\left\{x_i=\frac{i}{n+1}\, : \, i \in \mathbb{Z}\right\}$;
	\item $w(x_i,x_j)=\begin{cases}
	1 & \mbox{if } |i-j|=1,\\
	0 & \mbox{otherwise}.
	\end{cases}$
\end{itemize}
Fix $n \in \N$, define $X^{(n)} :=X\cap (0,1)$ and set the Dirichlet transformed-path graph Laplacian $\Delta_{\infty,\textnormal{dir}}$, associated to $X^{(n)}$ and $G$, with 
\begin{equation*}
\phi(m):= (-1)^{m+1} \frac{2}{m^2}.
\end{equation*}
The choice of $\phi$ is not random. Indeed, $\sum_{m\geq1}\phi(k)= - \pi^2/3$ and $\{\pi^2/3\}\cup \{-\phi(m)\}_{m\in\N}$ are the Fourier coefficients of $\zeta(\theta)=\theta^2$ over the interval $[0,\pi]$, and the inverse Weyl distribution function of the operator $\mathcal{L}$ in \eqref{eq:diff_operator} is exactly given by $\zeta$. Now, let us observe that $d_\infty=\infty$ and that $\Delta_{m,\textnormal{dir}} = 2I_n$ for every $m >n$, where $I_n$ is the identity matrix $n\times n$. It is not difficult to prove then that $\Delta_{\infty,\textnormal{dir}}$ is well-defined, and in particular let us observe that it is a symmetric Toeplitz matrix $n\times n$ whose stencil $\boldsymbol{t}$ is given by 
\begin{equation}\label{eq:Fourier_coeff}
\boldsymbol{t}= \left[\frac{\pi^2}{3},-2,\ldots,(-1)^n\frac{2}{n^2}\right].
\end{equation}  
This last remark will play a crucial role, as we will see shortly. For a review about Toeplitz matrices and their spectra, we refer to \cite{BG}. Finally, let us define the operator $\mathcal{L}^{(n)}$ via
\begin{equation*}
\mathcal{L}^{(n)}:= n^2\Delta_{\infty,\textnormal{dir}} + Q^{(n)},
\end{equation*}
where $Q^{(n)}$ is the $n\times n$ diagonal matrix whose entries are given by the pointwise evaluation of $q$ over $X^{(n)}$. It is a fact that $\mathcal{L}^{(n)}$ approximates $\mathcal{L}$ over $C_c^\infty(0,1)$, and it can be checked by two different approaches. Indeed, the stencil $\boldsymbol{t}$ can be obtained by both the sinc collocation method, see \cite{LK}, and as a limit of the $(2p+1)$-points Finite Difference method, where $p\to \infty$, see \cite{Bianchi}. If we now write $\delta_m^{(n)}$ and $\delta_m$ for the $m$-th eigenvalue of $\mathcal{L}^{(n)}$ and $\mathcal{L}$, respectively, then we obtain that
$$
\limsup_{n\to \infty} \max_{m=1,\ldots,n}\left\{\left|\frac{\delta^{(n)}_m}{\delta_m} -1\right|\right\}=0,
$$
that is, $\MSRE(\mathcal{L})=0$. This is due to this specific discretization and the stencil $\boldsymbol{t}$, since $\zeta(\theta)=\theta^2$ is exactly the (inverse) Weyl distribution function associated to $\mathcal{L}$, as we observed earlier. We do not get into the details here, for some references see \cite[Section 3.2]{Bianchi}. 

We can now define
\begin{equation}\label{eq:graph_approx}
\Kd := \mathcal{L}^{(n),\dagger},
\end{equation}
see Algorithm \ref{algo:1}. 
\begin{algorithm}\caption{Construction of $\Kd$ defined in \eqref{eq:graph_approx} by applying \eqref{eq:Green_operator}-\eqref{eq:diff_operator}}
	\begin{algo}\label{algo:1}
		\INPUT $n \in \N$, $q: (0,1)\to \mathbb{R}$
		\OUTPUT $\Kd \in \mathbb{R}^{n\times n}$
		\STATE Construct $\boldsymbol{q}:=\left[q\left(\frac{1}{n+1}\right),\ldots,q\left(\frac{n}{n+1}\right)\right]$
		\STATE Construct $Q^{(n)}= \textnormal{diag}(\boldsymbol{q}) \in \mathbb{R}^{n\times n}$
		\STATE Construct $\Delta_{\infty,\textnormal{dir}}= \textnormal{toeplitz}(\boldsymbol{t})\in \mathbb{R}^{n\times n}$, $\boldsymbol{t}$ as in \eqref{eq:Fourier_coeff}
		\STATE Define $\mathcal{L}^{(n)}= n^2\Delta_{\infty,\textnormal{dir}}+Q^{(n)}$
		\STATE Compute $U,\Sigma,V$ such that $\mathcal{L}^{(n)}= U\Sigma V^T$
		\STATE $\Kd = V\Sigma^\dagger U^T$
	\end{algo}
\end{algorithm}
The operator $\Kd$ such defined guarantees both \ref{item:local_convergence} and \ref{item:uniform_eig_estimate}. We can check it numerically through an example. Let $\K$ in \eqref{eq:Green_operator} be characterized by 
\begin{equation*}
h(x,y)= \begin{cases}
\sin^{-1}(1)\sin(1-x)\sin(y) & \mbox{if } 0\leq y < x \leq 1,\\
\sin^{-1}(1)\sin(x)\sin(1-y) & \mbox{if } 0\leq x\leq y\leq 1,
\end{cases}
\end{equation*}
or, equivalently, let $\mathcal{L}$ in \eqref{eq:diff_operator} be characterized by $q(x)\equiv -1$. We want to compare the discretization $\Kd$ obtained by \eqref{eq:graph_approx} with the discretization $\hat{\K}^{(n)}$ obtained by the Galerkin method with orthonormal box functions from \eqref{eq:Green_operator}. In Table \ref{table:Max_abs_error} and Table \ref{table:LSRE} it is possible to check the validity of \ref{item:local_convergence} and \ref{item:local_eig_estimate}, respectively, for both the discretizations. 
\begin{table}[H]
	\caption{\begin{flushleft}\textnormal{In this table we validate numerically the property \ref{item:local_convergence} for both the discretizations $\Kd$ and $\hat{\K}^{(n)}$. We choose the sup norm and $f(x)=x$ as test function. As it can be checked, the maximum absolute error tends to zero as $n$ increases, for both the discretizations.}\end{flushleft}}\label{table:Max_abs_error}
	\centering
	\begin{tabular}{cccc}
		\hline
		\multicolumn{4}{c}{Maximum Absolute Error}                                                           \\ 
		\hline
		\multicolumn{1}{l}{}      & $n=10^2$  & $n=10^3$ & $n=2\cdot10^3$  \\ 
		\hline
		$\mathcal{K}^{(n)}$       &   4.8094e-04      &      4.8213e-05           &      2.4110e-05          \\ 
		\hline
		$\hat{\mathcal{K}}^{(n)}$ &      0.0016             &    1.7749e-04    &      8.9109e-05      \\
		\hline
	\end{tabular}
\end{table}

\begin{figure}
	\centering
	\includegraphics[width=9cm]{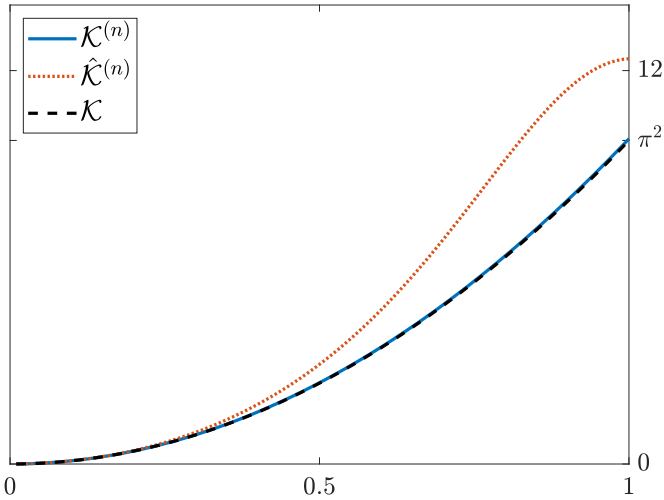}
	\caption{Graphical representation of (a suitable modification of) the first $n=100$ eigenvalues associated to $\K$, $\hat{\K}^{(n)}$ and $\Kd$. We have plotted the reciprocal of their eigenvalues, normalized by $n^2$ and sorted in nondecreasing order. The continuous blue line and the continuous red line are associated to the modified eigenvalues of $\hat{\K}^{(n)}$ and $\Kd$, respectively. The dashed black line is associated to the modified eigenvalues of $\K$. On the $x$-axis is reported the ratio $m/n$, where $m$ is referred to the $m$-th eigenvalue. As it can be seen, there is perfect match between the plots correspondent to $\K$ and $\Kd$, while the plot correspondent to $\hat{\K}^{(n)}$ deviates consistently from $\K$ in the mid-high frequencies. This behaviour confirms the results reported in Table \ref{table:MSRE}.}\label{fig:eig_comparison}
\end{figure}

\begin{table}[H]
\caption{\begin{flushleft}\textnormal{In this table we report the LSRE for both the discretizations $\Kd$ and $\hat{\K}^{(n)}$. As it can be checked numerically, the LSRE tends to zero as $n$ increases for fixed $m$, for both the discretizations, therefore validating property \ref{item:local_eig_estimate}.}\end{flushleft}}
\label{table:LSRE}
\centering
\begin{tabular}{ccccc} 
	\hline
	\multicolumn{5}{c}{LSRE}                                                                                     \\ 
	\hline
	&        & $n=10^2$ & $n=10^3$ & $n=2\cdot10^3$  \\ 
	\hline
	\multirow{3}{*}{$\mathcal{K}^{(n)}$ }      & $m=1$  &     0.0053                  &     5.3341e-04     &          2.6686e-04       \\ 
	& $m=10$  &    0.0048                   &     4.7988e-04    &      2.4007e-04        \\ 
	& $m=50$ &       0.0048                 &     4.7950e-04     &            2.3985e-04     \\ 
	\hline
	\multirow{3}{*}{$\hat{\mathcal{K}}^{(n)}$} & $m=1$  &  0.0100                   &    9.9983e-04      &          4.9996e-04       \\ 
	& $m=10$  &      0.0183                 &       0.0011   &               5.2034e-04  \\ 
	& $m=50$ &     0.1884                  &      0.0031    &              0.0010   \\
	\hline
\end{tabular}
\end{table}
The difference between the two discretizations relies in Property \ref{item:uniform_eig_estimate}. The eigenvalues of the continuous operator $\K$ are given by $\lambda_m = \frac{1}{m^2\pi^2-1}$, and Table \ref{table:MSRE} shows that $\Kd$ satisfies \ref{item:uniform_eig_estimate} while $\hat{\K}^{(n)}$ fails it. 
\begin{table}[H]
\caption{\begin{flushleft}\textnormal{In this table we report the MSRE defined in \eqref{def:MSRE}, for the two discretizations $\Kd$ and $\hat{\K}^{(n)}$. As it can be checked numerically, the MSRE associated to $\Kd$ tends to zero as $n$ increases while the MSRE associated to $\hat{\K}^{(n)}$ is stuck away from zero.}\end{flushleft}}\label{table:MSRE}
\centering
	\begin{tabular}{ccccc}
		\hline
		\multicolumn{5}{c}{MSRE}                                                           \\ 
		\hline
		\multicolumn{1}{l}{}      & $n=10^2$ & $n=5\cdot10^2$ & $n=10^3$ & $n=2\cdot10^3$  \\ 
		\hline
		$\mathcal{K}^{(n)}$       &     0.0053     &        0.0011        &    5.3341e-04      &       2.6686e-04          \\ 
		\hline
		$\hat{\mathcal{K}}^{(n)}$ &     0.3200     &        0.3098       &    0.3085      &         0.3078        \\
		\hline
	\end{tabular}
\end{table}
To give a visual understanding of what happens, in Figure \ref{fig:eig_comparison} we compare the eigenvalue distributions of $\Kd$ and $\hat{\K}^{(n)}$ with the eigenvalue distribution of $\K$. For a clear representation, we plotted the reciprocal of the first $n$ eigenvalues normalized by $n^2$.

\section{Graph Laplacian and the penalty term}\label{sec:graph_penalty}
We build the operator $\mathcal{A}^{(n)}$, in the penalty term of \eqref{def:discrete_generalized_Tikhonov}, such that it can represent the graph Laplacian associated to a specific graph $G$ that encodes informations from the observed data $g^\epsilon$. A very common choice is the following, see for example \cite{Peyre,Calatroni}: having fixed $n\in \N$, $r\in \{1,\ldots,n\}$ and $\sigma>0$, define $X=\left\{x_i=\frac{i}{n+1}\,:\, i=1,\ldots,n\right\}$ and 
\begin{equation}\label{def:w_penalty_term}
w(x_i,x_j) = \begin{cases}
\textnormal{e}^{-\frac{(g^\epsilon(x_i)-g^\epsilon(x_j))^2}{\sigma^2}} & \mbox{if } |i-j|\leq r,\\
0 &\mbox{otherwise}.
\end{cases}
\end{equation}
Then we write $\mathcal{A}^{(n)}:= \Delta$, where $\Delta$ is the graph Laplacian \eqref{formal_laplacian} associated to $G=(X,w)$ just defined above, see Algorithm \ref{algo:2}. The edge-weight function expresses the similarities between nodes, which (in this case) is given by a Gaussian distribution-like on the image-set of $g^\epsilon$. This choice revealed to be fruitful in a recent work on image deblurring, see \cite{BBDR}. Indeed, the kernel of the graph Laplacian $\Delta$, when the potential term $\kappa$ is zero, is one-dimensional and its eigenspace is generated by the constant vector $[1,\ldots,1]$. Images typically are characterized by having wide portions of its area constituted by almost homogeneous color, and then it happens that
\begin{equation*}
\left\|\Delta\left[f^\dagger\right]\right\| \approx 0,
\end{equation*}      
that is, $f^\dagger$ ``approximately'' belongs to the kernel of $\Delta$. Anyway, even if this approach generally can produce good results, there is still room for a great margin of improvement. Future lines of research have to concentrate on different $w$ that can take into account other players, such as the operator $\K$ itself, and a-priori informations regarding the true solution.  

\begin{algorithm}\caption{Construction of $\mathcal{A}^{(n)}$ in \eqref{def:discrete_generalized_Tikhonov}}
	\begin{algo}\label{algo:2}
		\INPUT $n \in \N$, $r\in \{1,\ldots,n\}$, $\sigma>0$, $\bsg_n \in \mathbb{R}^n$
		\OUTPUT $\mathcal{A}^{(n)} \in \mathbb{R}^{n\times n}$
		\STATE Construct $W\in \mathbb{R}^{n\times n}$ s.t. $(W)_{i,j}= w(x_i,x_j)$ as in \eqref{def:w_penalty_term}
		\STATE Construct $\boldsymbol{d}_n\in \mathbb{R}^n$ s.t.  $(\boldsymbol{d}_n)_i= \sum_{j=1}^n w(x_i,x_j)$
		\STATE Construct $D=\textnormal{diag}(\boldsymbol{d}_n)\in \mathbb{R}^{n\times n}$
		\STATE Define $\Delta= D-W$
		\STATE $\mathcal{A}^{(n)} =\Delta$
	\end{algo}
\end{algorithm}

\section{Numerical examples}\label{sec:numerical_examples}
In this section we provide some numerical examples of reconstructed solutions of the model problem \eqref{model_problem} by means of the discretized version of the regularization operator \eqref{def:gen_tik1}-\eqref{def:gen_tik2}, that is, 
$$
\bsf_{\alpha,n}^\epsilon:= \R_\alpha^{(n)} [\bsg_n^{\epsilon}]
$$
where $\R_\alpha^{(n)}$ is defined by \eqref{def:discrete_generalized_Tikhonov}. The continuous operator $\K$ will be an integral operator of the form \eqref{eq:Green_operator}, and we are going to consider three model examples in the upcoming subsections. For all the examples, the matrix operator $\mathcal{A}^{(n)}$ in the penalty term will be chosen among the identity matrix $I$ and the following ones
\begin{equation*}
\mathcal{A}_1^{(n)}= \begin{bmatrix}
2 & \shortminus 1 & 0 & \cdots & 0 \\
\shortminus 1 & 2 & \shortminus 1 &  & \\
  &   \ddots		&  \ddots   &       \ddots    &          \\
    &   		&     &           &          \\
  0    &   	\cdots	&   0  &       \shortminus 1    &         2 
\end{bmatrix},
\end{equation*}
\begin{equation*} 
\mathcal{A}_2^{(n)}= \begin{bmatrix}
1 & \shortminus 1 & 0 & \cdots & 0 \\
\shortminus 1 & 2 & \shortminus 1 &  & \\
&   \ddots		&  \ddots   &       \ddots    &          \\
&   		&     &           &          \\
0    &   	\cdots	&   0  &       \shortminus 1    &         1 
\end{bmatrix},
\end{equation*}
and $\mathcal{A}_3^{(n)}$, given by Algorithm \ref{algo:2}. Matrices $\mathcal{A}_1^{(n)}$ and $\mathcal{A}_2^{(n)}$ are the (normalized) discretization, by means of the $3$-points Finite Difference method, of the one-dimensional Laplace operator on $(0,1)$ with Dirichlet and Neumann BCs, respectively. The main operator $\Kd$ will be given by Algorithm \ref{algo:1}, and its performance will be compared with a standard discretization by the Galerkin method with orthonormal box functions. We point out that the latter discretization do not satisfy \ref{item:uniform_eig_estimate}.

All the computations are performed on \textsc{Matlab} R2020b. The noisy data vectors $\bsg_n^\epsilon$ are obtained by adding Gaussian-noise to the original data vectors $\bsg_n$ and, in order to have repeatability of the numerical experiments, we fixed $\texttt{rng(7)}$ for the function \texttt{randn()}. Specifically,
\begin{equation*}
\bsg_n^\epsilon = \bsg_n + \epsilon \frac{\mbox{\texttt{norm}(}\bsg_n \mbox{\texttt{)}}\mbox{\texttt{randn(}}1,n\mbox{\texttt{)}}}{\mbox{\texttt{norm(randn(}}1,n\mbox{\texttt{))}}},
\end{equation*}  
where $\epsilon$ is the noise level and $n$ is the number of points that (uniformly) discretize the interval $(0,1)$. The regularized solution $\bsf_{\alpha,n}^\epsilon$ is then obtained by applying the \texttt{tikhonov()} function which can be found in the \texttt{regtools} toolbox, see \cite{Hansen}, and the goodness of the reconstruction is evaluated by computing the \emph{relative restoration error} (RRE) in the $\ell^2$-norm,
\begin{equation*}
\textnormal{RRE}:= \frac{\|\bsf_{\alpha,n}^\epsilon - \bsf^\dagger_n\|}{\| \bsf^\dagger_n\|}.
\end{equation*}
In this work we do not focus on the strategies for choosing the regularization parameters $\alpha, r, \sigma$. In particular, $r$ and $\sigma$, that appears in \eqref{def:w_penalty_term}, can be difficult to set and ideally they should also be adapted to the noise level $\epsilon$. For this reason, for all the numerical examples we fix $r=\lceil(0.2)n\rceil$ and $\sigma=0.01$. The same goes for the $\alpha$ parameter: in all the experiments we choose the $\alpha$ that minimizes the RRE among fifty logarithmically spaced points between $10^3$ and $10^{-6}$.

\subsection{Example 1}\label{ssec:ex1}
In this first example we consider as Green kernel the following function
\begin{equation*}
h(x,y)= \begin{cases}
y(x-1) &\mbox{if } 0\leq y<x\leq 1,\\
x(y-1) &\mbox{if } 0\leq x\leq y\leq 1,
\end{cases}
\end{equation*}
which is the Green function of the operator \eqref{eq:diff_operator}, changed of sign and with $q(x) \equiv 0$. It is taken from the model problem \texttt{deriv2()} in the \texttt{regtools} toolbox. We consider the following test functions:
 \begin{equation*}
 f^\dagger_1(x) = \begin{cases}
 [p_2^2(x)-p_3(x)]\textrm{e}^{4-\frac{1}{p_1(x)}} & \mbox{if } (x-1/2)^2\leq 1/4,\\
 0  & \mbox{otherwise},
 \end{cases}
  \end{equation*}
  where
  \begin{equation*}
  p_1(x):= 0.25-(x-0.5)^2, \quad p_2(x):= 2\frac{(x-0.5)}{p_1^2(x)},\\
  \end{equation*}
  \begin{equation*}
   p_3(x):= \frac{2p_1^2(x)^2 + 8(x-0.5)^2p_1(x)}{p_1^4(x)};
  \end{equation*}
and
\begin{equation*}
f^\dagger_2(x)= \frac{x^3}{3}  - \frac{x^2}{2}, \quad f^\dagger_3(x)= x,\quad f^\dagger_4(x) = \textrm{e}^x.
\end{equation*}
Some preliminary remarks: we choose $f^\dagger_1$ because it belongs to $C^\infty_c(0,1)$, that is, the core subset of the operator \eqref{eq:diff_operator}, and $f^\dagger_2$ because it satisfies Neumann BCs. Instead, $f^\dagger_3$ and $f^\dagger_4$ are test functions originally implemented in  \texttt{deriv2()}. 

In the first two tables, we present two unrealistic examples but that are enlightening. They provide an extreme and clear confirmation of the statements we did in the previous Sections \ref{sec:graph_approx}-\ref{sec:graph_penalty}. In Table \ref{tab:epsilon=0}, we choose $f^\dagger_1$ and fixed $\epsilon=0$: since $f^\dagger_1 \in C^\infty_c(0,1)$, as we expected, a discretization of $\K$ that better preserve all the spectrum improves dramatically the RRE. In Table \ref{tab:true_solution} instead, we fixed $\epsilon=0.1$ and assumed to know already the true solution $f^\dagger_3$, and we used this information to build a potential term $\kappa$ to add to the graph Laplacian $\mathcal{A}^{(n)}_3$ such that $\mathcal{A}^{(n)}_3\left[\bsf^\dagger_{n,3}\right]=0$. Despite the high level of noise, the RRE is very low, as it can be compared with the results in Table \ref{tab:example1}, where a lower lever of noise is implemented. This is due to the fact that the true solution lives in the kernel of $\mathcal{A}^{(n)}_3$. It is interesting to notice how strong is the projection of the regularized solution $\bsf_{\alpha,n}^\epsilon$ into the kernel of $\mathcal{A}^{(n)}_3$, so that the role of the discretization of $\K$ becomes secondary.

\begin{table}[!htb]
	\begin{minipage}{.5\linewidth}
		\caption{\begin{flushleft}\textnormal{Toy-model example with $f^\dagger_1$ and $\mathcal{A}^{(n)}=I$. The noise level is $\epsilon=0$ and $n=100$.} \end{flushleft}}\label{tab:epsilon=0}
		\centering
		\begin{tabular}{cc} 
			\hline
			\multicolumn{2}{c}{RRE}  \\ 
			\hline
			$\mathcal{K}^{(n)}$       &             3.9195e-07      \\ 
			\hline
			$\hat{\mathcal{K}}^{(n)}$ &             0.0187      \\
			\hline
		\end{tabular}
	\end{minipage}%
	\begin{minipage}{.5\linewidth}
		\centering
		\caption{\begin{flushleft}\textnormal{Toy-model example with $f^\dagger_3$ and $\mathcal{A}^{(n)}$ such that $\mathcal{A}^{(n)}[\bsf^\dagger_{n,3}]=0$. The noise level is $\epsilon=0.1$ and $n=100$.} \end{flushleft}}\label{tab:true_solution}
	\begin{tabular}{cc} 
		\hline
		\multicolumn{2}{c}{RRE}  \\ 
		\hline
		$\mathcal{K}^{(n)}$       &             0.0021      \\ 
		\hline
		$\hat{\mathcal{K}}^{(n)}$ &             0.0024      \\
		\hline
	\end{tabular}
	\end{minipage} 
\end{table}
In Table \ref{tab:example1} we provide the RRE for several combinations of $\Kd, \hat{\K}^{(n)}$ and $\mathcal{A}^{(n)}$. The first remark is that for every fixed operator $\mathcal{A}^{(n)}$ in the penalty term, the best RRE is given by using $\Kd$, for every $f^\dagger_i$. The second remark is that the best RRE is gained by the pair $\{\Kd, \mathcal{A}^{(n)}_3\}$, for all the test cases with the only exception of $f^\dagger_1$. In the test cases $f^\dagger_i$, for $i=2,3,4$, the improvement is remarkable while in the test case $f^\dagger_1$ the presence of $\mathcal{A}^{(n)}_3$ worsen drastically the RRE. This is most probably due to the fact that the choice made for $w$ in \eqref{def:w_penalty_term} does not fit well with functions that are compactly supported in $(0,1)$.

\begin{table}[H]
	\centering
	\caption{\begin{flushleft}\textnormal{For all these numerical experiments we fixed $\epsilon=0.01$ and $n=100$. In bold are highlighted the best RRE for each $f^\dagger_i$, $i=1,2,3,4$.} \end{flushleft}}
	\label{tab:example1}
	\begin{tabular}{cccccl} 
		\hline
		\multicolumn{6}{c}{RRE}                                                                                                                                       \\ 
		\hline
		\multicolumn{1}{l}{}           & \multicolumn{1}{l}{}      & $I$ & $\mathcal{A}^{(n)}_1$ & $\mathcal{A}^{(n)}_2$ & \multicolumn{1}{c}{$\mathcal{A}^{(n)}_3$}  \\ 
		\hline
		\multirow{2}{*}{$f^\dagger_1$} & $\mathcal{K}^{(n)}$       &  0.0836   &       \textbf{0.0664}                &        0.0869               &               0.2430                             \\
		& $\hat{\mathcal{K}}^{(n)}$ &   0.0876  &            0.0686           &        0.0930               &                   0.2416                          \\ 
		\hline
		\multirow{2}{*}{$f^\dagger_2$} & $\mathcal{K}^{(n)}$       &   0.2184  &             0.2275        &       0.0185                &       \textbf{0.0176}                                    \\
		& $\hat{\mathcal{K}}^{(n)}$ &    0.2262 &                 0.2477        &             0.0308           &       0.0313                                      \\ 
		\hline
		\multirow{2}{*}{$f^\dagger_3$} & $\mathcal{K}^{(n)}$       &   0.2017  &            0.1954           &     0.0391                  &                 \textbf{0.0328}                           \\
		& $\hat{\mathcal{K}}^{(n)}$ &   0.2098  &           0.2258            &      0.0772                 &               0.0678                             \\ 
		\hline
		\multirow{2}{*}{$f^\dagger_4$} & $\mathcal{K}^{(n)}$       & 0.1937    &      0.1912                 &   0.0293                    &           \textbf{0.0089}                                 \\
		& $\hat{\mathcal{K}}^{(n)}$ &   0.2021  &            0.2227           &        0.0453               &       0.0150                                     \\
		\hline
	\end{tabular}
\end{table}

\begin{figure}[H]
	\centering
	\includegraphics[width=9cm]{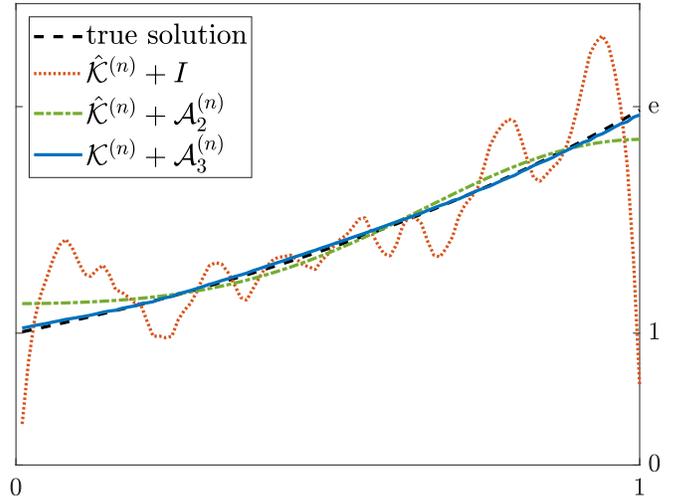}
	\caption{Plots of the regularized solutions $\bsf_{\alpha,n}^\epsilon$ from the test case $f^\dagger_4$ in Table \ref{tab:example1}, for several different combinations. The dashed-black line represents the true solution $f^\dagger_4$. The best RRE is achieved by the pair $\Kd$ and $\mathcal{A}^{(n)}_3$, that is represented by the continuous-blue line, and that can be easily checked by a direct visual inspection. Observe that if we choose $\mathcal{A}^{(n)}=I$ then $\bsf_{\alpha,n}^\epsilon$ is forced to assume zero values at the endpoints $x=0,x=1$. This is due to the fact that $\K$ is the Green operator associated to $\mathcal{L}$ defined in \eqref{eq:diff_operator}, and $\textnormal{dom}\left(\mathcal{L}\right)=H^1_0= \overline{C^\infty_c(0,1)}$.}\label{fig:ex1}
\end{figure}

Finally, in Figure \ref{fig:ex1} we report the plots of the regularized solutions $\bsf_{\alpha,n}^\epsilon$ from the test case $f^\dagger_4$ in Table \ref{tab:example1}, for several different combinations. 

\subsection{Example 2}
\begin{figure}[H]
	\centering
	\includegraphics[width=9cm]{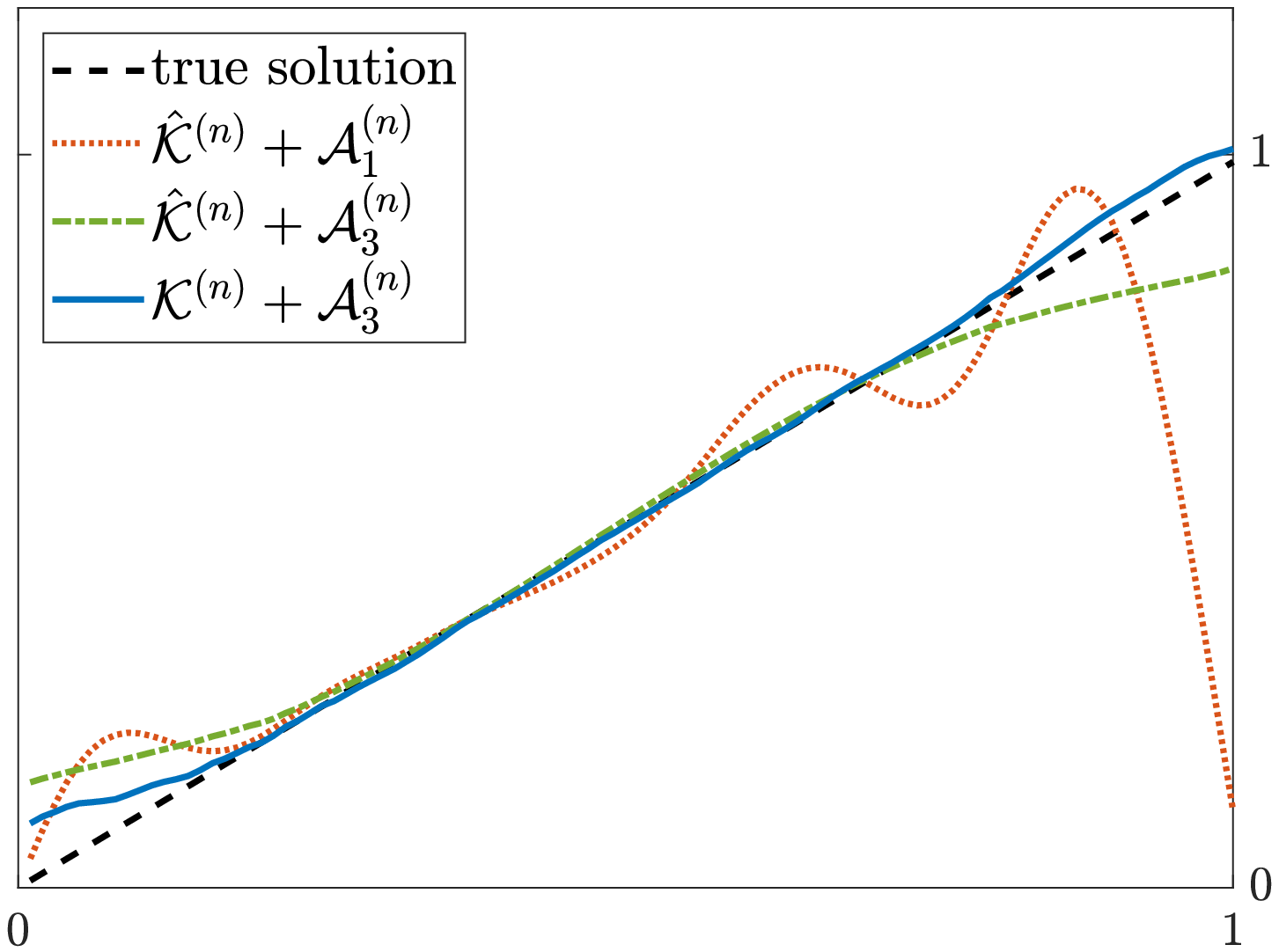}
	\caption{Plots of the regularized solutions $\bsf_{\alpha,n}^\epsilon$ from the test case $f^\dagger_3$ in Table \ref{tab:example2}, for several different combinations. The dashed-black line represents the true solution $f^\dagger_3$. The best RRE is achieved by the pair $\Kd$ and $\mathcal{A}^{(n)}_3$, that is represented by the continuous-blue line, and that can be easily checked by a direct visual inspection. Observe that if we choose $\mathcal{A}^{(n)}_1$ then $\bsf_{\alpha,n}^\epsilon$ is forced to assume zero values at the endpoints $x=0,x=1$. This is due to the fact that $\mathcal{A}^{(n)}_1$ is a (normalized) discretization of $\mathcal{L}$ defined in \eqref{eq:diff_operator} for $q=0$, and $\textnormal{dom}\left(\mathcal{L}\right)=H^1_0= \overline{C^\infty_c(0,1)}$.}\label{fig:ex2}
\end{figure}

In this second example we take
\begin{equation*}
h(x,y)= \begin{cases}
\sin^{-1}(1)\sin(1-x)\sin(y) & \mbox{if } 0\leq y < x \leq 1,\\
\sin^{-1}(1)\sin(x)\sin(1-y) & \mbox{if } 0\leq x\leq y\leq 1,
\end{cases}
\end{equation*}
as in Section \ref{sec:graph_approx}, and we keep the $f^\dagger_i$, for $i=1,\ldots,4$ defined in the preceding Subsection \ref{ssec:ex1}, as test functions. We collected the numerical results in Table \ref{tab:example2}, while in Figure \ref{fig:ex2} we provided the plots of the reconstructed solutions for the test case $f^\dagger_3$, for several different combinations of the discretized operator and the penalty term.
\begin{table}[H]
	\centering
	\caption{\begin{flushleft}\textnormal{For all these numerical experiments we fixed $\epsilon=0.02$ and $n=100$. In bold are highlighted the best RRE for each $f^\dagger_i$, $i=1,2,3,4$.} \end{flushleft}}
	\label{tab:example2}
	\begin{tabular}{cccccl} 
		\hline
		\multicolumn{6}{c}{RRE}                                                                                                                                       \\ 
		\hline
		\multicolumn{1}{l}{}           & \multicolumn{1}{l}{}      & $I$ & $\mathcal{A}^{(n)}_1$ & $\mathcal{A}^{(n)}_2$ & \multicolumn{1}{c}{$\mathcal{A}^{(n)}_3$}  \\ 
		\hline
		\multirow{2}{*}{$f^\dagger_1$} & $\mathcal{K}^{(n)}$       & 0.1202    &   \textbf{0.0893}                   &            0.1133          &      0.4557                                    \\
		& $\hat{\mathcal{K}}^{(n)}$ &    0.1190  &         0.0905            &    0.1138               &       0.4636                                      \\ 
		\hline
		\multirow{2}{*}{$f^\dagger_2$} & $\mathcal{K}^{(n)}$       &   0.2404  &     0.2185              &               0.0169         &         \textbf{0.0161}                                 \\
		& $\hat{\mathcal{K}}^{(n)}$ &   0.3330  &     0.3306                &     0.0636                 &              0.0623                             \\ 
		\hline
		\multirow{2}{*}{$f^\dagger_3$} & $\mathcal{K}^{(n)}$       &  0.2637   &   0.2553                    &         0.0426             &      \textbf{0.0406}                                    \\
		& $\hat{\mathcal{K}}^{(n)}$ &   0.3200  &          0.2908           &    0.1075                   &     0.1003                                       \\ 
		\hline
		\multirow{2}{*}{$f^\dagger_4$} & $\mathcal{K}^{(n)}$       &  0.2455   &   0.2395                    &        0.0307              &           \textbf{0.0116}                              \\
		& $\hat{\mathcal{K}}^{(n)}$ &   0.2937 &         0.2703            &    0.0900                 &              0.0712                             \\
		\hline
	\end{tabular}
\end{table}

\section{Conclusions}
The main purpose of this short paper is to show that a good combination, of a discretization technique for the main operator $\K$ and the operator $\mathcal{A}$ in the penalty term, can improve by far the reconstructed solutions, if compared to some standard discretizations. We used a graph-based approach to both build $\Kd$ and $\mathcal{A}^{(n)}$, even if the two discretizations had different scopes. In the first case, we were interested to have a good spectral relative error all along the spectrum, in comparison with the spectrum of the real operator $\K$, while in the second case we wanted to find a good space where to force into the reconstructed solution. The preliminary numerical results are promising, nevertheless there are still many open questions. For example, it is necessary to find a procedure to build a discretization of $\K$ from the kernel $h$ that guarantees $\MSRE=0$, without knowing the differential operator $\mathcal{L}$. Moreover, it has to be understood if all the spectrum should be preserved or if it is possible to concentrate on a smaller portion. Finally, the edge-weight function $w$, used to build the graph Laplacian $\mathcal{A}^{(n)}$ in the penalty term, should be refined such to encode more possible informations on the space where the true solution $f^\dagger$ lives.


\begin{thebibliography}{99}
\bibitem{AH} H. C. Andrews, B. R. Hunt, {\em Digital image restoration}. Prentice–Hall, Englewood Cliffs, NJ (1977). 
\bibitem{ABS} A. Adriani, D. Bianchi, S. Serra-Capizzano, {\em Asymptotic Spectra of Large (Grid) Graphs with a Uniform Local Structure (Part I): Theory}. Milan J. Math. 88(2) (2020): 409--454.
\bibitem{Bianchi} D. Bianchi, {\em Analysis of the spectral symbol associated to discretization schemes of linear self-adjoint differential operators}. Preprint (2020), arXiv:2004.10058.  
\bibitem{BD} D. Bianchi, M. Donatelli, {\em On generalized iterated Tikhonov regularization with operator-dependent seminorms}. ETNA 47 (2017): 73--99.
\bibitem{BBDR} D. Bianchi, A. Buccini, M. Donatelli, E. Randazzo, {\em Graph Laplacian for image deblurring}. Preprint (2021), arXiv:2102.10327.
\bibitem{Davies1} E. B.	Davies, \emph{Spectral theory and differential operators}. Cambridge University Press (1996).
\bibitem{BG} A. B\"{o}ttcher, S. M. Grudsky, {\em Toeplitz matrices, asymptotic linear algebra and functional analysis}. Vol. 67. Springer (2000).
\bibitem{Calatroni} L. Calatroni, Y. van Gennip, C. B. Sch{\"{o}}nlieb, H. M. Rowland, A. Flenner, {\em Graph Clustering, Variational Image Segmentation Methods and Hough Transform Scale Detection for Object Measurement in Images}. J. Math. Imaging Vis. 57(2) (2017): 269--291.
\bibitem{EHN} H. W. Engl, M. Hanke, A. Neubauer, {\em Regularization of inverse problems}. Vol. 375. Springer Science \& Business Media (1996).
\bibitem{Estrada1} E. Estrada, {\em Path Laplacian matrices: Introduction and application to the analysis of consensus in networks}. Linear Algebra Appl. 436(9) (2012): 3373--3391.
\bibitem{Estrada2} E. Estrada, E. Hameed, N. Hatano, M. Langer, {\em Path Laplacian operators and superdiffusive processes on graphs. I. One-dimensional case}. Linear Algebra Appl. 523 (2017): 307--334.
\bibitem{EM} W. N. Everitt, L. Markus, \emph{Boundary value problems and symplectic algebra for ordinary differential
	and quasi-differential operators}. American Mathematical Soc. (1999).
\bibitem{G} C. W. Groetsch, {\em The theory of tikhonov regularization for fredholm equations}. Boston Pitman Publication (1984).
\bibitem{Hansen} P. C. Hansen, {\em  Regularization tools: a Matlab package for analysis and solution of discrete illposed problems}. Numer. Algorithms 6 (1994): 1--35.
\bibitem{HS} T. K. Huckle, M. Sedlacek, {\em Tikhonov-Phillips regularization with operator dependent seminorms}. Numer. Algorithms 60 (2012): 339--353.
\bibitem{KL} M. Keller, D. Lenz, {\em Dirichlet forms and stochastic completeness of graphs and subgraphs}. J. fur die Reine und Angew. Math. 666 (2012): 189--223.
\bibitem{KLW_book} M. Keller, D. Lenz, R. K. Wojciechowski, \emph{Graphs and discrete Dirichlet spaces}. Grundlehrender mathematischen Wissenschaften: Springer, forthcoming (2021).
\bibitem{LNEM} F. Lenti, F. Nunziata, C. Estatico, M. Migliaccio, {\em Spatial resolution enhancement of earth observation products using an acceleration technique for iterative methods}. IEEE Geosci. Remote Sensing 12 (2015): 269--273.

\bibitem{LK} J. Lund, K. L. Bowers, {\em Sinc methods for quadrature and differential equations}. Society for Industrial and Applied Mathematics (1992).
\bibitem{Peyre} G. Peyr{\'{e}}, S.  Bougleux, L. Cohen, {\em Non-local regularization of inverse problems}. In: Computer Vision -- ECCV 2008, Springer Berlin Heidelberg, Eds D. Forsyth, P. Torr, A. Zisserman (2008): 57--68.
\end{thebibliography}
\end{document}